\documentclass[9pt]{article}
\usepackage{spconf,amsmath,graphicx}
\usepackage[hashEnumerators,smartEllipses]{markdown}
\usepackage{soul}
\usepackage{xcolor}
\usepackage[ruled]{algorithm2e}
\usepackage{cleveref}
\usepackage{caption}
\usepackage{subcaption}

\def\R{{\mathbb R}}
\DeclareMathOperator*{\argmin}{arg\,min}

\newtheorem{definition}{Definition}
\title{Robust Data-Driven Accelerated Mirror Descent}
\name{\begin{tabular}{c}Hong Ye Tan$^{\star}$ \qquad Subhadip Mukherjee$^{\star \dagger}$ \qquad Junqi Tang$^{\star}$ \\
Andreas Hauptmann$^{\ddagger \S}$ \qquad Carola-Bibiane Sch\"onlieb$^{\star}$\end{tabular}}

\address{$^{\star}$ University of Cambridge, DAMTP, Cambridge CB3 0WA, U.K. \\
$^{\dagger}$ Department of Computer Science, University of Bath, U.K. \\
$^{\ddagger}$ University of Oulu, Research Unit of Mathematical Sciences, P.O.Box 8000, 90014 University of Oulu \\
$^{\S}$ University College London, Department of Mathematics, 25 Gordon St, London WC1H 0AY, U.K.
}
\begin{document}
%
\maketitle
\begin{abstract}

Learning-to-optimize is an emerging framework that leverages training data to speed up the solution of certain optimization problems. One such approach is based on the classical mirror descent algorithm, where the mirror map is modelled using input-convex neural networks. In this work, we extend this functional parameterization approach by introducing momentum into the iterations, based on the classical accelerated mirror descent. Our approach combines short-time accelerated convergence with stable long-time behavior. We empirically demonstrate additional robustness with respect to multiple parameters on denoising and deconvolution experiments.
\end{abstract}
%


    


\section{Introduction}
\label{sec:intro}

Learning-to-optimize is a general framework of methods that seek to minimize some objective functions quickly (in as few iterations as possible). This has been proposed in various settings such as unsupervised learning or using primal-dual methods, with applications such as imaging inverse problems and neural network training \cite{li2016learningtooptimize,banert2020, chen2021L2O}. A recent work focuses on introducing a learned functional parameterization into the classical mirror descent (MD) algorithm, obtaining faster convergence rates with approximate convergence guarantees \cite{LMD22Tan}.

In this work, we propose a learned convex optimization approach based on the classical accelerated mirror descent (AMD) scheme. Our aim is to minimize a family of convex functions $f$ in some function class $\mathcal{F}$ of qualitatively similar problems, such as variational image denoising. In this work, we will study and demonstrate convergence and robustness of our proposed learned scheme, and compare it to existing learned and classical optimization methods.  \Cref{sec:AMD} outlines the AMD scheme, and \cref{sec:robustness} contains various experimental results on the robustness of our proposed method. 


\subsection{Mirror Descent}

Mirror descent is a generalization of gradient descent (GD), first introduced by Nemirovsky and Yudin \cite{nemirovskij1983problem}. By exploiting the geometry of the cost function, MD achieves competitive convergence rate bounds for certain problems, including online learning and tomographic reconstruction \cite{orabona2015generalized, tomographyMD2001bental}. We first outline the method as presented by Beck and Teboulle \cite{BECK2003167}.

Let $\Psi$ be a convex function on a closed convex set $\mathcal{X} \subseteq \R^n$. Let $(\R^n)^*$ denote the corresponding dual space of $\R^n$. We denote by $\bar{\R} = \R \cup \{+\infty\}$ the extended real line. Recall that the \emph{convex conjugate} (or \emph{Fenchel conjugate}) of $\Psi$, denoted $\Psi^* : (\R^n)^* \rightarrow \bar{\R}$, is given by $\Psi^*(y) = \sup_{x \in \mathcal{X}} (\langle y, x \rangle - f(x))$. 
$\Psi$ induces a \emph{Bregman divergence} $B_\Psi : \mathcal{X} \times \mathcal{X} \rightarrow \bar{\R}$, defined as:
    \[B_\Psi(x,y) = \Psi(x) - \Psi(y) - \langle \nabla \Psi(y), x-y \rangle.\]
    
\begin{definition}[Mirror potential]
    We say $\Psi:\mathcal{X} \rightarrow \R$ is a \textbf{mirror potential} if it is continuously differentiable and strongly convex. We call the gradient $\nabla \Psi:\mathcal{X} \rightarrow (\R^n)^*$ a \textbf{mirror map}.
\end{definition}
If $\Psi$ is a mirror potential, then the convex conjugate $\Psi^*$ is everywhere differentiable, and moreover satisfies $\nabla \Psi^* = (\nabla \Psi)^{-1}$ \cite{BECK2003167, RockWets1998VA}. The standard MD update rule for minimizing convex $f$ with initialization $x_0 \in \mathcal{X}$ is as follows \cite{BECK2003167}:
\begin{equation}\label{eq:MDDef}
    \begin{split}
        y_{k} = \nabla \Psi (x_k) - t_k \nabla f (x_k), \,\,\, x_{k+1} = \nabla \Psi^* (y_k).
    \end{split}
\end{equation}
For convex Lipschitz $f$, MD is able to attain $\mathcal{O}(1/\sqrt{k})$ convergence rate. Similarly to GD, algorithmic extensions such as stochasticity, ergodicity or acceleration are available for MD \cite{accelMD15Walid, accelmd2018hanzely,ergodicMD2012duchi, stochasticMD2021orazio}. MD can be additionally be interpreted as natural GD, with a Riemannian metric induced by the mirror map \cite{mdgeom2015garvesh, mirrorlessmd2021gunasekar}. While the inclusion of the mirror map makes this method more flexible than GD, current applications are limited by the requirement of a hand-crafted mirror map. For efficient computation, such a mirror map is further restricted by the requirement of a closed-form convex conjugate, which is generally unavailable. 

\subsection{Learned MD}\label{sec:LMD}
One way to make the MD algorithm learnable and alleviate the requirement of a hand-crafted mirror map is to replace $\nabla \Psi$ and $\nabla \Psi^*$ with learnable maps $\nabla M_\theta$ and $\nabla M_\vartheta^*$. The learned mirror potential $M_\theta$ is parameterized as an input-convex neural network (ICNN) to enforce convexity \cite{amos2017icnn}. While a closed-form convex conjugate is generally unavailable for an ICNN $M_\theta$, this method learns another function $M_{\vartheta}^*$ to approximate the convex conjugate.

To maintain a MD structure for this scheme, the learned mirror maps are enforced to be close inverses of each other, $\nabla M_\vartheta^* \circ \nabla M_\theta \approx I$. We will refer to the distance between $\nabla M_\vartheta^* \circ \nabla M_\theta$ and the identity $I$ as the \emph{forward-backward error}. Using such a parameterization results in the learned MD (LMD) method with step-sizes $(t_k)_{k\ge 0}$ detailed in \cite{LMD22Tan}:
\begin{equation}
    \tilde{x}_{k+1} = \nabla M_\vartheta^* (\nabla M_\theta (\tilde{x}_k) - t_k \nabla f(\tilde{x}_k)).
\end{equation}
The goal of LMD is to accelerate convergence on a class of functions, typically with qualitatively similar attributes, such as image denoising or inpainting. LMD has been shown to outperform methods such as GD and Adam for various convex problem classes. Formally, fix a function class $\mathcal{F}$ consisting of convex functions $f:\mathcal{X} \rightarrow \R$, where $\mathcal{X}\subseteq \R^d$. The goal is to minimize functions in this class $f \in \mathcal{F}$ quickly on average via LMD, with initializations $\tilde{x}_0 = x$ distributed possibly depending on $f$. For training, one possible loss function is to minimize the function values up to the $N$-th iteration for some fixed $N$, while also penalizing the distance of the mirror maps from being inverses of each other. This enforces the resulting LMD method to optimize quickly on the function class, and stay close to a proper MD method for convergence guarantees. This can be expressed as a minimization problem over the mirror potential parameters $\theta, \vartheta$:
\begin{equation}\label{eq:MDTrainLoss}
    \underset{\theta,\vartheta}{\min}\,\mathbb{E}_{f,x}\left[\sum_{k=1}^N r_k f(\tilde{x}_k) + s_k\|(\nabla M_\vartheta^* \circ \nabla M_\theta - I)(\tilde{x}_k)\|\right].
\end{equation}

\section{Accelerated MD}\label{sec:AMD}
The faster convergence rates of the Nesterov accelerated GD method can be extended to MD \cite{Nesterov1983AMF}. This can be done by considering the dynamics of certain ODEs corresponding to MD and Nesterov accelerated methods, and combining them in a suitable manner to derive an accelerated MD ODE \cite{nemirovskij1983problem,accelMD15Walid, differential2014su}. Discretizing the resulting ODE and applying a small modification results in a family of accelerated MD (AMD) methods, which we summarize in the following section.

Let $\mathcal{X} \subseteq \R^n$ be a closed convex set, and $f:\mathcal{X} \rightarrow \bar{\R}$ be a proper $C^1$ convex function. Suppose further that $\nabla f$ is $L_f$-Lipschitz (w.r.t. a norm $\|\cdot\|$), and let $f^*$ be the minimum of $f$ on $\mathcal{X}$. Assume that $\Psi^*$ is $L_{\Psi^*}$-smooth with respect to $\|\cdot\|_*$, the dual norm of $\|\cdot\|$ on the dual space. Let $R$ be a regularization function such that there exists $0 < \ell_R \le L_R$ such that for all $x, x' \in \mathcal{X}$, $\frac{\ell_R}{2}\|x-x'\|^2 \le R(x, x') \le \frac{L_R}{2}\|x-x'\|^2$. The resulting AMD algorithm with initialization $\tilde{x}^{(0)} = x_0, \tilde{z}^{(0)} = \nabla \Psi^* (x_0)$ is as follows:
\begin{subequations}
\begin{align}\label{eq:AMD}
        &\lambda_k = \frac{r}{r+k},\\
        &x^{(k+1)} = \lambda_k\tilde{z}^{(k)} + (1-\lambda_k) \tilde{x}^{(k)}, \\
        &\tilde{z}^{(k+1)} = \nabla \Psi^*(\nabla\Psi(\tilde{z}^{(k+1)}) - \frac{kt}{r} \nabla f(x^{(k+1)})), \\
        &\tilde{x}^{(k+1)} = \argmin_{\tilde{x} \in \mathcal{X}} \gamma t \langle \nabla f(x^{(k+1)}), \tilde{x} \rangle + R(\tilde{x}, x^{(k+1)}). \label{eq:amdgdstep}
\end{align}
\end{subequations}
\subsection{Learned AMD}
We propose to make the AMD method learnable by replacing the mirror maps $\nabla \Psi$ and $\nabla \Psi^*$ with learned mirror maps $\nabla M_\theta$ and $\nabla M_\vartheta^*$, respectively, as explained in \Cref{sec:LMD}. We additionally learn the step-sizes $t_k$ for the first $N=10$ iterations, similar to the adaptive LMD training in \cite{LMD22Tan}. For simplicity, we fix $r=3$ and $\gamma=1$. We take the regularizer $R(x, x') = \frac{1}{2}\|x-x'\|_2^2$, turning (\ref{eq:amdgdstep}) into a GD step. The resulting learned iterates are as in \Cref{alg:learnedAMD}.

\begin{algorithm}
\DontPrintSemicolon
\caption{Learned AMD (LAMD)}\label{alg:learnedAMD}
\KwData{Mirror potential $\Psi$, step-sizes $(t_k)_{k=1}^N>0$, parameter $r\ge 3$}
Initialize $\tilde{x}^{(0)}= x_0, \tilde{z}^{(0)} = x_0$.\;
\For{$1\le k \le N$}{
    $\lambda_k = \frac{r}{r+k}$.\;
    $x^{(k+1)} = \lambda_k\tilde{z}^{(k)} + (1-\lambda_k) \tilde{x}^{(k)}$\;
    $\tilde{z}^{(k+1)} = \nabla M_\vartheta^*(\nabla M_\theta(\tilde{z}^{(k+1)}) - \frac{kt_k}{r} \nabla f(x^{(k+1)}))$\;
    $\tilde{x}^{(k+1)} = x^{(k+1)} - \gamma t_k \nabla f(x^{(k+1)})$\;
}%
\end{algorithm}\vspace*{-.4cm}
\subsection{Training procedure}
The first function class that we train on is for TV denoising of noisy ellipse phantoms in X-ray CT \cite{wang2006total,landi2012efficient}. The ellipse phantoms dataset was generated using the Deep Inversion Validation Library (DIVal) \cite{leuschner2019deep}. A sinogram is first created from an ellipse phantom of size $128^2$ by first applying bilinear upsampling to $400^2$ to avoid inverse crime \cite{kaipio2006statistical}, then applying a parallel-beam ray transform with 30 angles and 183 measurements per angle, and finally adding 10$\%$ Gaussian noise to the sinogram. A noisy phantom is created by taking the filtered back-projection (FBP) of this sinogram and downsampling it to $128^2$. Denoting the noisy phantom thus generated by $y$, the resulting function class of TV-regularized variational losses is:
\begin{equation}
    \mathcal{F} = \left\{f(x) = \|x-y\|_2^2 + \lambda\|\nabla x\|_{1}\  \right\},
\end{equation}
where $\lambda$ is a regularization parameter, the gradient is taken pixel-wise, and the norms are taken over the image space. This setup allows us to train on realistic noise artifacts without introducing an expensive forward operator into training. We chose $\lambda=0.15$ by visually comparing the TV-based reconstructions after running GD for 1000 iterations. The initializations were chosen to be the noisy phantoms $x_0 = y$.


For an initialization $x_0$, let $(x^{(k)})_{k=1}^N$ be the iterates generated by LAMD \Cref{alg:learnedAMD}. The loss function that we wish to minimize for training will be of the form (\ref{eq:MDTrainLoss}) trained for $N=10$ iterations, to minimize the function values at the iterates $f(x^{(k)})$, as well as to promote good mirror map inversion around the iterates and the ground-truths. 
\section{Robustness study of LAMD}\label{sec:robustness}
When using LMD or LAMD as a potential alternative for optimizing convex functions, the method should converge after running for many iterations. In the following section, we investigate the convergence behavior of the learned methods with different step-sizes, and generality with respect to the objective function class. 
\subsection{Extension of learned step-sizes}
We first extend the LMD and LAMD iterations past the learned $N=10$ iterations and investigate the long-time behavior of the iterations. This extension is done by continuing the corresponding MD update steps with the same learned mirror maps $M_\theta, M_\vartheta^*$, and by extending the step-sizes $t_k$ using various methods. The extension methods used are based on the learned step-sizes $t_1,...,t_{10}$. In particular, we use the maximum, mean, and minimum of the learned step-sizes, the final learned step-size $t_k = t_{10}$, and a ``reciprocal" step-size. The reciprocal step-size is $t_k = c/k$, which was observed to be close to the learned step-sizes. Here, $c$ was taken to fit the first $N=10$ iterations, i.e. $c = \frac{1}{N} \sum_{k=1}^{N} k t_k$.

\begin{figure}[t!b]
\centering
    \begin{picture}(230,120)
    \put(0,-10){\includegraphics[width=8.5cm]{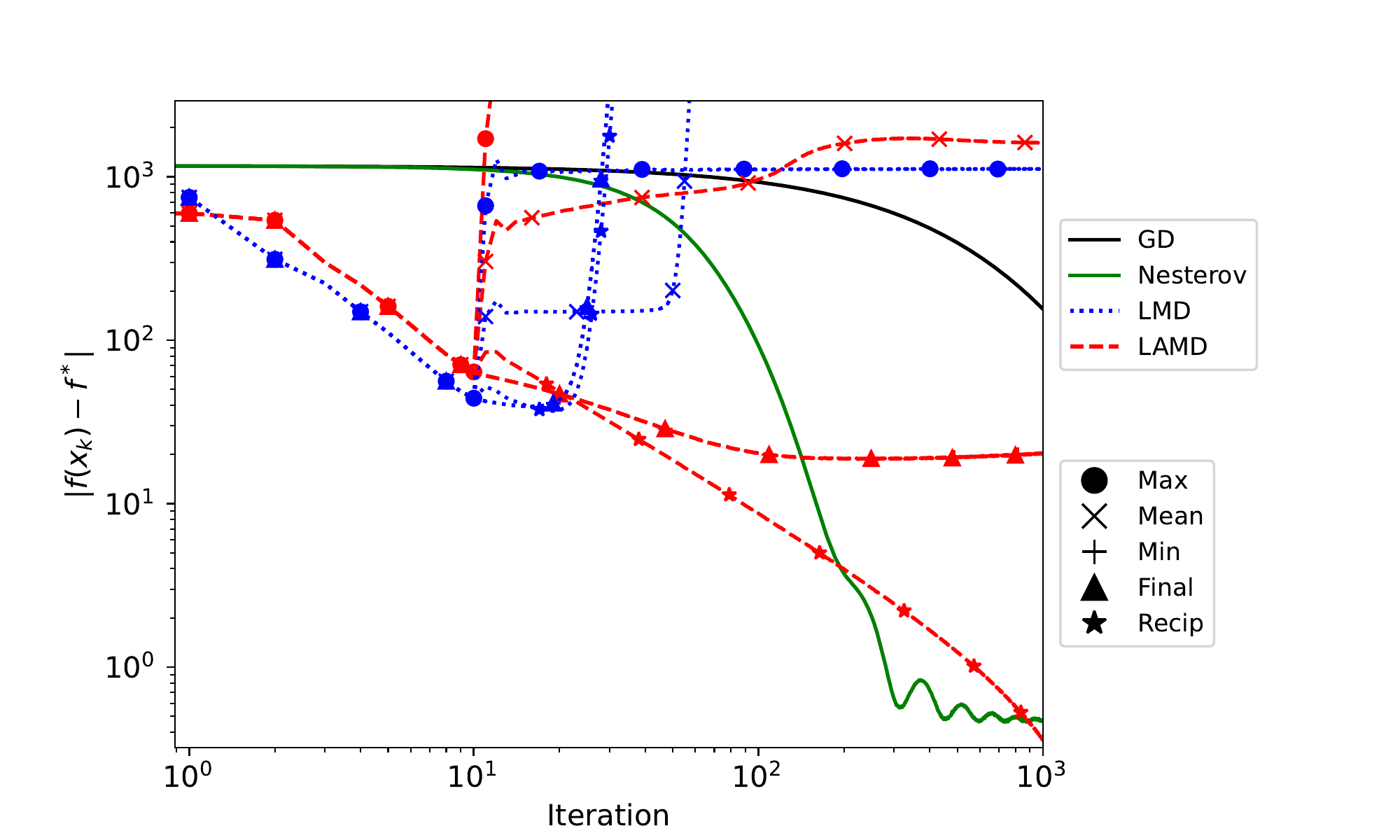}}
\end{picture}
\caption{Plot of $f(x^{(k)}) - f^*$ with various step-size extensions, where $f^* = \min f$. Extending step-sizes past the 10 learned iterations with the choice $t_k = c/k$ gives the best convergence. The minimum learned step-size here is the final one.}
\label{fig:stepsize_ext_logloss}
\end{figure}

We observe in \Cref{fig:stepsize_ext_logloss} that the non-accelerated LMD method suffers from instability when taking small step-size extensions, and poor performance for large step-size extensions. While approximate convergence guarantees are available in this setting, a necessary condition is that the forward-backward error is uniformly bounded on the iterates, which is enforced for the trained iterates. However, this condition is violated when going past the trained number of iterates, as there is no guarantee that the forward-backward error stays low. In contrast, LAMD does well when taking small step-size extensions, as shown by the reciprocal extensions attaining roughly $\mathcal{O}(1/k)$ convergence rate. For the reciprocal extension, the iterates appear to continue decreasing the loss long after the learned number of iterations, making this a clear choice of step-size for further experiments.


\vspace{-.4cm}

\subsection{Domain transfer}

We propose to transfer the LMD and LAMD methods away from denoising ellipse phantoms, introducing various structural changes to the objective function class. In particular, we will consider changing the ray transform, using a different phantom dataset, and using a convolution forward operator. 

Changing the ray transform changes the noise distribution on the phantoms, introducing different artifacts. In particular, the parallel-beam ray transform in the noisy phantom generation is changed to a cone-beam ray transform with 60 angles and 400 measurements per angle. The parallel- and cone-beam transform experiments were also run on the LoDoPaB dataset, a benchmark for low dose CT reconstruction \cite{lodopab2021leuschner}.

\begin{figure*}[ht!]
\begin{picture}(450,126)
  \centering
   \put(10,-8){\begin{subfigure}[b]{8.5cm}
   \centering
  \includegraphics[width=8.5cm]{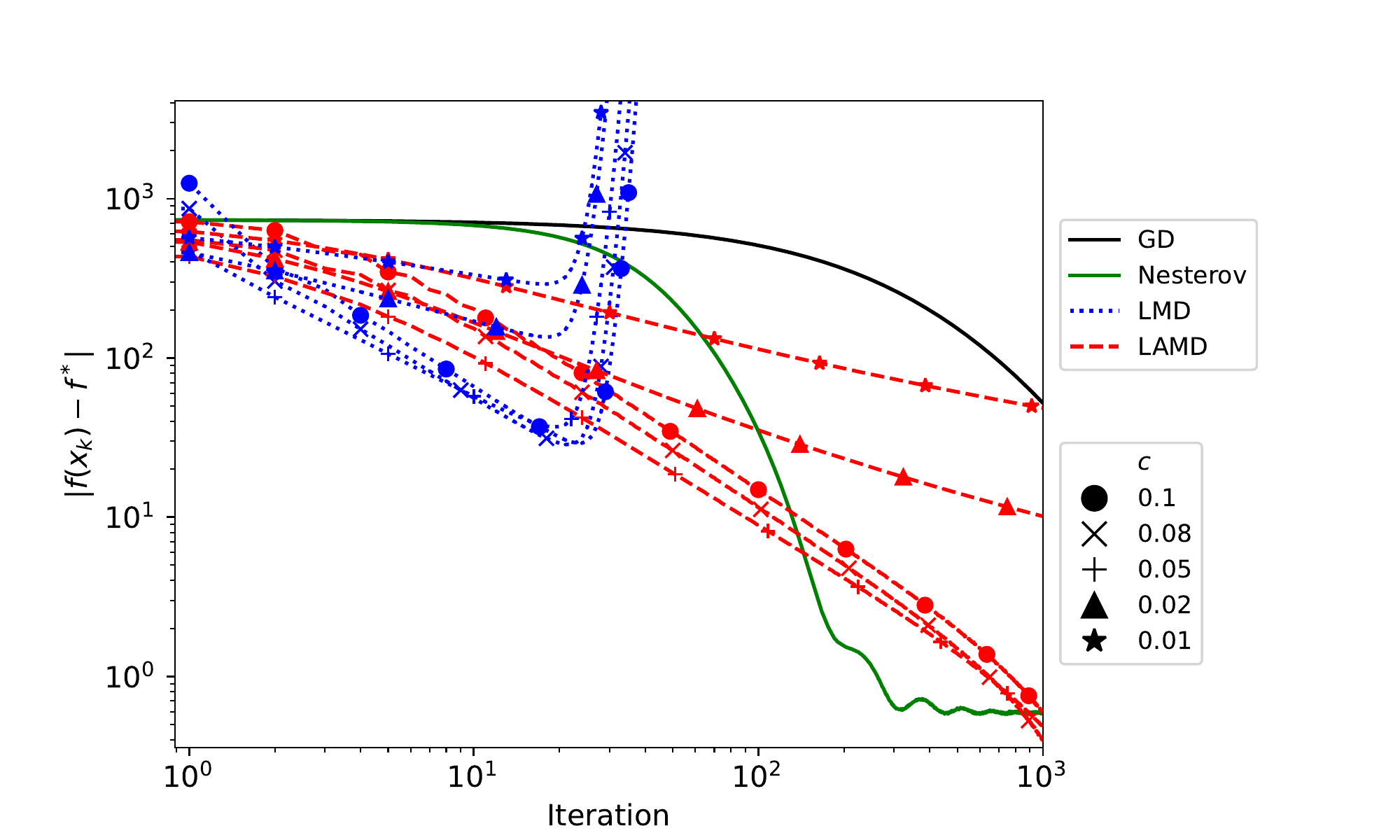}
  \vspace{-0.6cm}
  \caption{Ellipses dataset.}
    \label{fig:ellipse_cone}
    \end{subfigure}}
 \put(260,-8){\begin{subfigure}[b]{8.5cm}
   \centering
   \includegraphics[width=8.5cm]{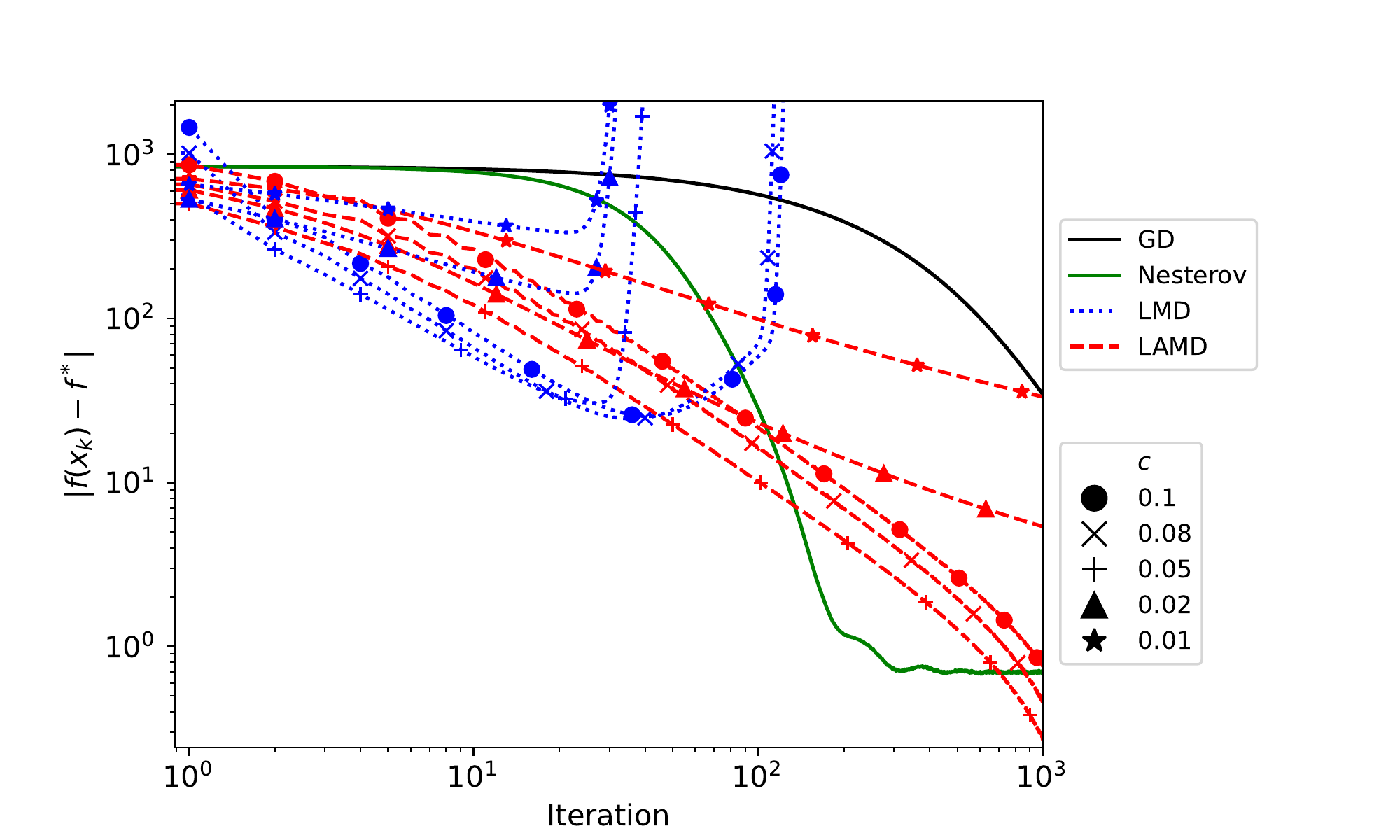}
    \vspace{-0.6cm}
    \caption{LoDoPaB dataset.}
    \label{fig:lodopab_res}
\end{subfigure}}
\end{picture}
\vspace{-0.1cm}
\caption{Plot of $f(x^{(k)})-f^*$ with various reciprocal step-sizes and cone beam transform. LMD and LAMD both generalize to denoising LoDoPaB phantoms, with LAMD achieving better convergence up to $10^3$ iterations. }
\end{figure*}

With the reciprocal step-sizes $t_k = c/k$, \Cref{fig:ellipse_cone,fig:lodopab_res} demonstrate that LMD continues to have good convergence for the earlier iterations, while still showing divergence as the iterations progress. LAMD continues to have roughly $\mathcal{O}(1/k)$ convergence rate in the later iterations, while still outperforming GD and the Nesterov accelerated scheme for the earlier iterations. We note that in Figure \ref{fig:lodopab_res}, LMD appears to diverge later for the LoDoPaB dataset than with the ellipses dataset. 

\begin{figure}[t!b]
\centering
    \begin{picture}(240,112)
    \put(0,-8){\includegraphics[clip, trim=0 0 0 40, width=8.5cm]{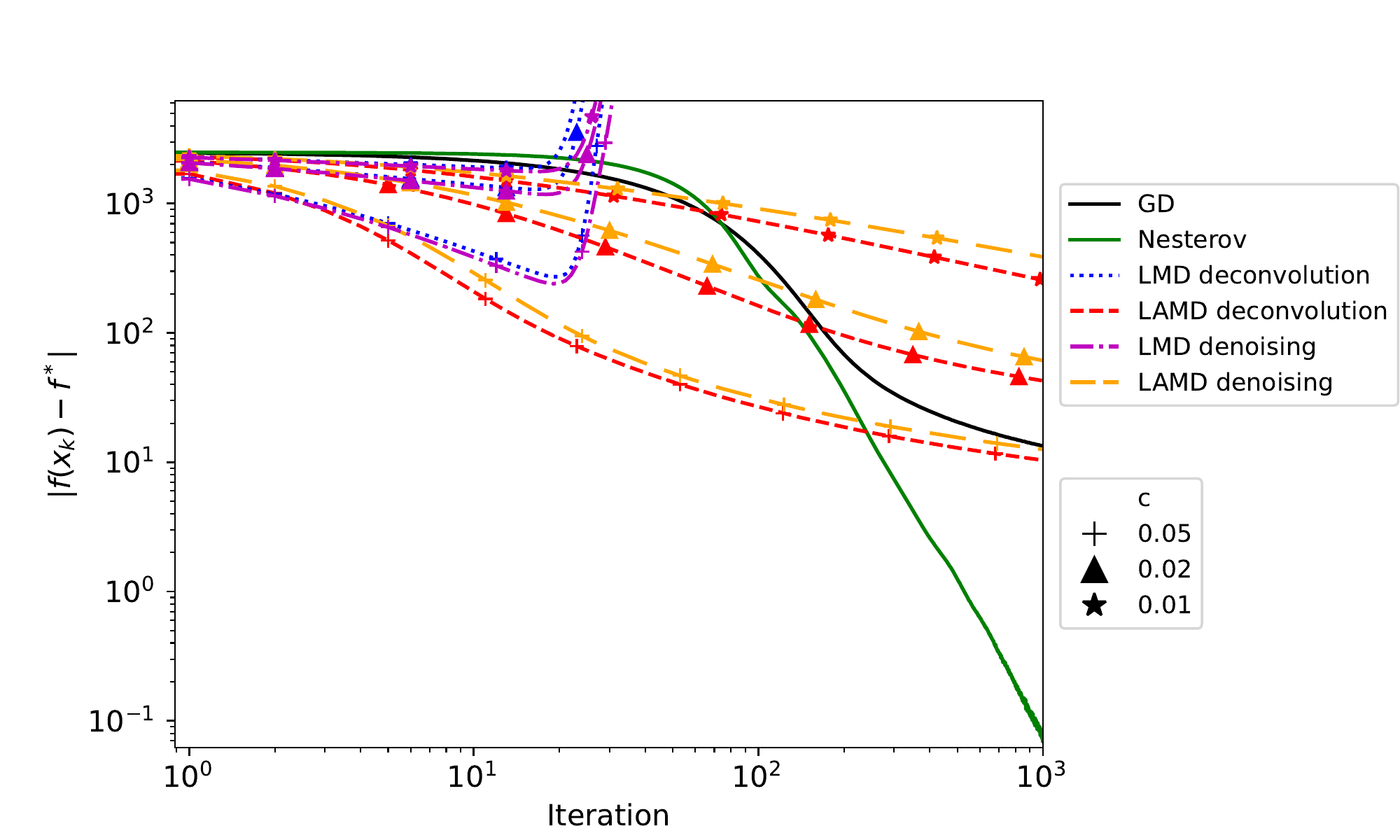}}
    \end{picture}
\caption{Plot of $f(x^{(k)})- f^*$ for deconvolution problem with various step-sizes, optimized using LMD/LAMD learned on deconvolution and on denoising. The similar performance of LMD and LAMD suggests a general learned mirror map which can be transferred across similar forward operators.}
\label{fig:ellipse_deconv}
\end{figure}

We also experiment with changing the forward operator. For the denoising experiments, the forward operator was taken to be the identity. Let $A$ be a Gaussian convolution kernel with kernel size of 7px and standard deviation of 3px. We generate blurred phantoms $y$ by applying the convolution $A$ and adding 10\% additive Gaussian noise. By changing the forward operator to the convolution $A$, we obtain a new function class corresponding to a deconvolution inverse problem:
\begin{equation}
    \mathcal{F}_{\text{conv}} = \left\{\|Ax - y\|^2 + \lambda \|\nabla x\|_1 : \text{blurred phantoms } y\right\}.
\end{equation}
\Cref{fig:ellipse_deconv} demonstrates the effect of applying LMD and LAMD with reciprocal step-sizes to the deconvolution function class, with mirror maps trained on either deconvolution or denoising problems. We see that LMD trained on deconvolution and denoising have very similar behavior. Moreover, for LAMD, the optimizer trained on denoising performs only marginally worse than that trained on deconvolution for the three choices of $c$ shown. This suggests generality of the mirror map, as applying the LAMD method trained on denoising to a different problem class still results in a stable and convergent method.




\subsection{Effect of momentum on training}
In this subsection, we investigate the difference between the mirror maps of the LMD and LAMD approaches. As the iterates are computed in a different manner, the training dynamics will also be different. Intuitively, we would expect the trained mirror maps to fit the geometry of the problem class. 



\begin{figure}[t!b]
\centering
    \begin{picture}(240,112)
    \put(0,-8){\includegraphics[clip, trim=0 0 0 40, width=8.5cm]{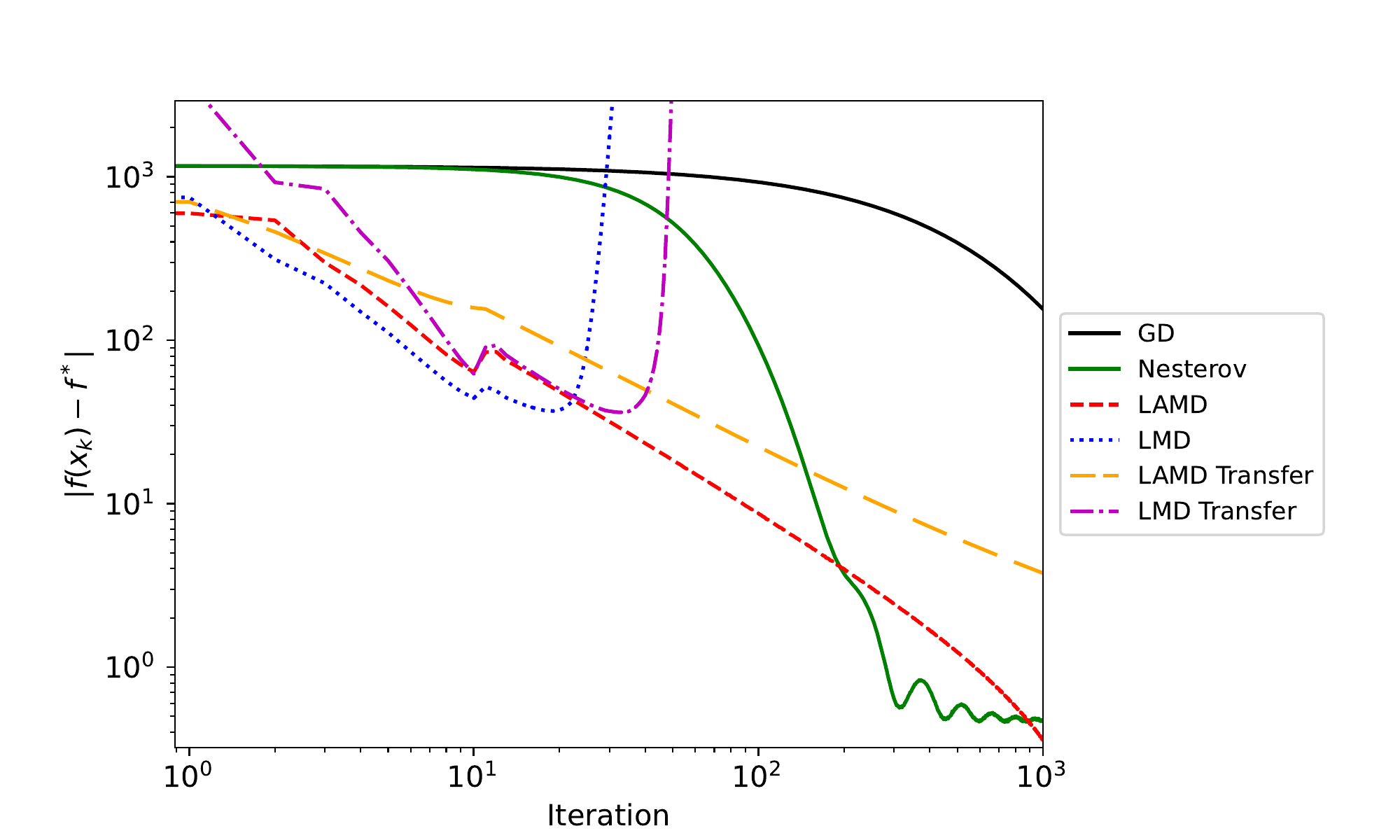}}
\end{picture}
\caption{Plot of $f(x^{(k)})-f^*$ for LMD and LAMD for denoising with swapped mirror maps. Adding momentum and using the mirror map of LMD (yellow) gives long-term stability. Training with LAMD also adds stability to the mirror maps when applied without momentum (purple).} %
\label{fig:ellipse_swap}
\end{figure}

We empirically compare the mirror maps by transferring the learned mirror maps between the LMD and LAMD methods. By transferring the mirror map learned using LMD to the LAMD scheme, or equivalently adding momentum to the LMD, we get the ``LAMD Transfer" scheme. Transferring the maps learned using LAMD to the LMD scheme yields the ``LMD Transfer" scheme. In \Cref{fig:ellipse_swap}, we observe that adding acceleration to the LMD method still results in reasonable convergence rates, albeit worse than the LAMD method. If momentum is removed from the LAMD method, the resulting iteration again performs worse than the trained LMD. However, the map trained using LAMD is more stable, diverging after some later iterations. This suggests that using momentum in the training process encourages the resulting mirror maps to be more stable, even with different training dynamics. Moreover, including momentum in the testing phase is a simple but effective method for adding stability to the iterates. 


\section{Conclusions}
In this work, we propose Learned Accelerated Mirror Descent, which outperforms methods such as GD and Nesterov accelerated GD \cite{LMD22Tan}. The LAMD scheme empirically results in stable extensions to further step-sizes, as well as more consistent mirror maps. The experiments also suggest that mirror maps trained on certain forward operators generalize to similar forward operators, such as from the identity to a convolution operator. A future avenue of research would be to extend convergence results for AMD to the approximate case where the inverse mirror map is not known exactly, or investigate robustness with respect to the momentum parameter. 



%
\vfill\pagebreak

\bibliographystyle{IEEEbib}
\bibliography{refs}

\end{document}